\documentclass[12pt]{elsart}
\usepackage{amsfonts,amsmath}

\journal{
}

\date{}

\begin{document}

\begin{frontmatter}
\title
{Type ${\rm III_1}$ factors generated by regular representations
of infinite dimensional nilpotent group $B_0^{\mathbb N}$}

\thanks[dec]{The  author would like to thank the Max-Planck-Institute
of Mathematics, Bonn
for the hospitality. The partial financial support by the DFG
project 436 UKR 113/87 is gratefully acknowledged. }

\author{Alexandre Kosyak
}
\ead{kosyak01@yahoo.com, kosyak@imath.kiev.ua}
\address{Max-Planck-Institut f\"ur Mathematik, Vivatsgasse 7, D-53111 Bonn,
Germany}
\address{Institute of Mathematics, Ukrainian National Academy of Sciences,
3 Tereshchenkivs'ka, Kyiv, 01601, Ukraine,\\
E-mail: kosyak01@yahoo.com, kosyak@imath.kiev.ua \\
tel.: 38044 2346153 (office), 38044 5656758 (home), fax: 38044
2352010}

\corauth[cor]{Corresponding author}

\newpage

\begin{abstract}
We study the von Neumann algebra, generated by the unitary
representations of infinite-dimensional groups nilpotent group
$B_0^{\mathbb N}$. The conditions of the irreducibility of the
regular and quasiregular representations of infinite-dimensional
groups (associated with some quasi-invariant measures) are given
by the so-called Ismagilov conjecture (see
\cite{AlbKos05C,AlbKos06,Kos90,Kos92,Kos94}). In this case the
corresponding von Neumann algebra is type ${\rm I}_\infty$ factor.
When the regular representation is reducible we find the
sufficient conditions on the measure for the von Neumann algebra
to be factor (see \cite{KosZek00,KosZek01}). In the present
article we determine the type of corresponding factors. Namely we
prove that the von Neumann algebra generated by the regular
representations of infinite-dimensional nilpotent group
$B_0^{\mathbb N}$ is type ${\rm III}_1$ hyperfinite factor. The
case of the nilpotent group $B_0^{\mathbb Z}$ of infinite in both
directions matrices will be studied in \cite{KosDyn07}.
\end{abstract}


\end{frontmatter}

\sloppy

\newcommand{\tr}{\mathrm{tr}\,}
\newcommand{\rank}{\mathrm{rank}\,}
\newcommand{\diag}[1]{\mathrm{diag}\,(#1)}
\renewcommand{\Im}{\mathrm{Im}\,}

\maketitle
\newpage
\tableofcontents



\section{ Regular representations} Let us consider the group $\tilde G=B^{\mathbb N}$ of all
upper-triangular real matrices of infinite order with unities on
the diagonal
$$
\tilde G=B^{\mathbb N}=\{I+x\mid x=\sum_{1\leq k<n}x_{kn}E_{kn}\},
$$
and its subgroup
$$
G=B_0^{\mathbb N}=\{I+x\in B^{\mathbb N}\mid\,\,x\,\,{\rm
is}\,\,{\rm finite}\},
$$
where $E_{kn}$ is an infinite-dimensional matrix with $1$ at the
place $k,n\in {\mathbb N}$ and zeros elsewhere,
 $x=(x_{kn})_{k<n}$ is {\it finite} means that $x_{kn}=0$ for
all $(k,n)$ except for a finite number of indices $k,n\in{\mathbb
N}$.

Obviously, $B_0^{\mathbb N}=\varinjlim_{n} B(n,{\mathbb R})$ is
the inductive limit of the group $B(n,{\mathbb R})$ of real
upper-triangular matrices with units on the principal diagonal
$$
B(n,{\mathbb R})=\{I+\sum_{1\leq k<r\leq n }x_{kr}E_{kr}\mid
x_{kr}\in {\mathbb R}\}
$$
with respect to the imbedding $B(n,{\mathbb R})\ni x\mapsto
x+E_{n+1n+1}\in B(n+1,{\mathbb R})$.

We define the Gaussian measure $\mu_b$ on the group $B^{\mathbb
N}$ in the following way
\begin{equation}
\label{mu-b} d\mu_b(x) = \otimes_{1\leq
k<n}(b_{kn}/\pi)^{1/2}\exp(-b_{kn}x_{kn}^2 )dx_{kn}=
\otimes_{k<n}d\mu_{b_{kn}}(x_{kn}),
\end{equation}
where $b=(b_{kn})_{k<n}$ is some set of positive numbers.

Let us denote by $R$ and $L$ the right and the left action of the
group $ B^{\mathbb N}$ on itself: $R_s(t)=ts^{-1},\,\,\,L_s(t)=
st, \,\,s,t \in B^{\mathbb N}$ and by $ \Phi:B^{\mathbb N}\mapsto
B^{\mathbb N},\,\,\Phi(I+x):=(I+x)^{-1}$ the inverse mapping.  It
is known \cite{Kos90,Kos92} that
\begin{lem}
$\mu_b^{R_t}\sim \mu_b\,\,\forall t\in B_0^{\mathbb
N}$  for any set $b=(b_{kn})_{k<n}$.
\end{lem}
\begin{lem}
$\mu_b^{L_t}\sim \mu_b\,\,\forall t\in B_0^{\mathbb N}$ if and
only if $S^L_{kn}(b)<\infty,\,\,\forall k<n$, where
$$
S^{L}_{kn}(b)=\sum_{m=n+1}^\infty\frac{b_{km}}{b_{nm}}.
$$
\end{lem}

\begin{lem}
$\mu_b^{L_t}\perp \mu_b\,\,\forall t\in B_0^{\mathbb
N}\backslash\{e\}\Leftrightarrow S^{L}_{kn}(b)=\infty \,\,\,
\forall k<n.$
\end{lem}

\begin{lem}{\rm \cite{Kos00}}
If $E(b)=\sum_{k<n}S^{L}_{kn}(b)(b_{kn})^{-1}<\infty$, then
$\mu_b^\Phi\sim\mu_b.$
\end{lem}
\begin{lem}{\rm \cite{Kos00}}
The measure $\mu_b$ on $B^{\mathbb N}$ is  $B_0^{\mathbb N}$
ergodic with respect to the right action.
\end{lem}
Let  $\alpha:G \rightarrow {\rm Aut}(X)$ be a measurable action of
a group $G$ on the measurable space $X$. We recall that a measure
$\mu$ on the space $X$ is $G$-ergodic if
$f(\alpha_t(x))=f(x)\,\,\forall t\in G$ implies
$f(x)=const\,\,\mu$ a.e. for all functions $f\in L^1(X,\mu)$.
\begin{rem}{\rm \cite{KosZek00}}
If $\mu_b^\Phi\sim\mu_b$ then $\mu_b^{L_t}\sim \mu_b\,\,\forall
t\in B_0^{\mathbb N}.$
\end{rem}
\begin{pf}
This follows from the fact that the inversion $\Phi$ replace the
right and the left action: $R_t\circ\Phi=\Phi\circ L_t\,\,\forall
t \in B^{\mathbb N}$. Indeed, if we denote
$\mu^{f}(\cdot)=\mu(f^{-1}(\cdot))$ we have
$(\mu^{f})^g=\mu^{f\circ g}$. Hence
$$
\mu_b\sim\mu_b^{R_t}\sim(\mu_b^{R_t})^\Phi= \mu_b^{R_t\circ\Phi}=
\mu_b^{\Phi\circ L_t} =(\mu_b^\Phi)^{L_t}\sim\mu_b^{L_t}.
$$
\qed\end{pf} If $\mu_b^{R_t}\sim \mu_b$ and $\mu_b^{L_t}\sim
\mu_b\,\,\forall t\in B_0^{\mathbb N}$, one can define in a
natural way (see \cite{Kos90,Kos92}), an analogue of the right
$T^{R,b}$ and left $T^{L,b}$ representation of the group
$B_0^{\mathbb N}$ in Hilbert space $H_b=L_2(B^{\mathbb N},d\mu_b)$
$$
T^{R,b},\,\,T^{L,b}:B_0^{\mathbb N}\rightarrow
U(H_b=L_2(B^{\mathbb N},d\mu_b)),
$$
$$(T^{R,b}_t f)(x)=(d\mu_b(xt)/d\mu_b(x))^{1/2}f(xt),$$
$$(T^{L,b}_s f)(x)=(d\mu_b(s^{-1}x)/d\mu_b(x))^{1/2}f(s^{-1}x).
$$
\section{ Von Neuman  algebras generated by the regular representations}
Let  ${\mathfrak A}^{R,b}=(T^{R,b}_t \mid t \in B_0^{\mathbb
N})^{\prime\prime}$ (resp.  ${\mathfrak A}^{L,b}=(T^{L,b}_s \mid s
\in B_0^{\mathbb N})^{\prime\prime}$) be the von Neumann algebras
generated by the right $T^{R,b}$  (resp. the left $T^{L,b}$)
regular representation of the group $B_0^{\mathbb N}$.
\begin{thm}{\rm \cite{Kos00}} If $E(b)<\infty$ then
$\mu_b^\Phi\sim\mu_b$. In this case the left regular
representation is well defined and the commutation theorem holds:
\begin{equation}
\label{(R)'=(L)} ({\mathfrak A}^{R,b})'={\mathfrak A}^{L,b}.
\end{equation}
 Moreover, the operator $J_{\mu_b}$ given by
\begin{equation}
\label{J}
(J_{\mu_b}f)(x)=(d\mu_b(x^{-1})/d\mu_b(x))^{1/2}\overline{f(x^{-1})}
\end{equation}
is an intertwining operator:
$$
T^{L,b}_t=J_{\mu_b}T^{R,b}_tJ_{\mu_b},\,\,t\in B_0^{\mathbb
N}\,\,\text{\, and\,}\,\, J_{\mu_b}{\mathfrak
A}^{R,b}J_{\mu_b}={\mathfrak A}^{L,b}.$$
\end{thm}
If $\mu_b^{L_t}\perp \mu_b\,\,\forall t\in B_0^{\mathbb
N}\backslash\{e\}$ one can't define the left regular
representation of the group $B_0^{\mathbb N}$. Moreover the
following theorem holds
\begin{thm} The right regular representation $T^{R,b}:B_0^{\mathbb N}\mapsto
U(H_b)$ is irreducible if and only if $\mu_{b}^{L_{s}}\perp
\mu_{b}\,\, \forall s\in B_0^{\mathbb N}\backslash\{0\}.$
\end{thm}
\begin{cor}
The von Neumann algebra ${\mathfrak A}^{R,b}$ is a type $I_\infty$
factor if \par
$\mu_{b}^{L_{s}}\perp \mu_{b}\,\, \forall s\in B_0^{\mathbb
N}\backslash\{0\}.$
\end{cor}
Let us  assume now that $\mu_b^{L_t}\sim \mu_b\,\,\forall t\in
B_0^{\mathbb N}\backslash\{e\}$. In this case the right regular
representation and the left regular representation of the group
$B_0^{\mathbb N}$  are well defined.

In \cite{KosZek00} the condition were studied  {\it when the von
Neumann algebra ${\mathfrak A}^{R,b}$ is factor}, i.e.
$${\mathfrak A}^{R,b}\cap({\mathfrak A}^{R,b})^\prime= \{\lambda
 {\bf I}\vert\lambda\in {\mathbb C}^1\}.
$$
Since $T^{L,b}_t \in ({\mathfrak A}^{R,b})^\prime \,\, \forall
t\in B_0^{\mathbb N}$, we have  ${\mathfrak
A}^{L,b}\subset({\mathfrak A}^{R,b})^\prime$, hence
\begin{equation}
{\mathfrak A}^{R,b}\cap({\mathfrak A}^{R,b})^\prime \subset
({\mathfrak A}^{L,b})^\prime\cap({\mathfrak A}^{R,b})^\prime =
({\mathfrak A}^{R,b}\cup{\mathfrak A}^{L,b})^\prime .
\end{equation}
The last relation shows that ${\mathfrak A}^{R,b}$ is factor if
 the representation
$$
B_0^{\mathbb N}\times B_0^{\mathbb N}\ni (t,\,s) \rightarrow
T^{R,b}_t T^{L,b}_s\in U(H_b)
$$
is irreducible.

Let us denote by ${\mathfrak A}^{R,L,b}$ the the von Neumann
algebras generated by the right $T^{R,b}$ and the left $T^{L,b}$
regular representations of the group $B_0^{\mathbb N}$:
$$
{\mathfrak A}^{R,L,b}=(T^{R,b}_t,T^{L,b}_s \mid t,s\in
B_0^{\mathbb N})''=({\mathfrak A}^{R,b}\cup {\mathfrak
A}^{L,b})''.
$$
Let us denote
$$
S^{R,L}_{kn}(b)=\sum_{m=n+1}^\infty
\frac{b_{km}}{S^{L}_{nm}(b)},\,\,\, k<n.
$$
\begin{thm}{\rm \cite{KosZek00}} The representation
$$
B_0^{\mathbb N}\times B_0^{\mathbb N}\ni (t,\,s) \rightarrow
T^{R,b}_t T^{L,b}_s\in U(H_b)
$$
is irreducible if $S^{R,L}_{kn}(b)=\infty,\,\,\forall k<n$.
\end{thm}
\begin{cor}
The von Neumann algebra ${\mathfrak A}^{R,b}$ is factor if
$S^{R,L}_{kn}\!(b)\!=\!\infty\,\forall k\!<\!n.$
\end{cor}

\section{Type ${\rm III}_1$ factor}
Let us denote as before  $M={\mathfrak A}^{L,b}=(T^{L,b}_s \mid s
\in B_0^{\mathbb N})^{\prime\prime},\,$ ${\mathfrak
A}^{R,b}=(T^{R,b}_t \mid t \in B_0^{\mathbb N})^{\prime\prime}$.
\begin{thm}
\label{III-0} If $S^{R,L}_{kn}(b)=\infty,\,\,\forall k<n$ then the
 von Neumann algebra ${\mathfrak A}^{L,b}$ (and hence
${\mathfrak A}^{R,b}$) is ${\rm III}_1$ factor.
\end{thm}
\begin{pf} The proof is based on Lemma \ref{Sp(De)} and
\ref{M-phi}, we shall prove them later.

Using (\ref{J}) we conclude that the modular operator $\Delta$ is
defined as follows
\begin{equation}
\label{Delta} (\Delta f)(x)=(d\mu_b(x)/d\mu_b(x^{-1}))f(x).
\end{equation}
\begin{lem}
\label{Sp(De)} We have
$$
Sp\Delta = [0,\infty).
$$
\end{lem}
We have $Sp\Delta\phi =Sp\Delta = [0,\infty),$ where
$\phi(a)=(a{\bf 1},\,{\bf 1})_{H_b},\,\,a \in M= {\mathfrak
A}^{L,b}.$ The centralizer $M_\phi$ of $\phi$ is defined by the
equality
$$
M_\phi= \{a \in M\mid \sigma^\phi_t(a)\,\,\forall t\in {\mathbb
R}\}
$$
where $\sigma^\phi_t(a)=\Delta^{it}a\Delta^{-it}.$ For every
projection $e\not= 0,\,e\in M_\phi ,$ a faithful semifinite normal
weight $\phi_e$ on the reduced von Neumann algebra $eMe=\{a\in
M;\,\,ea=ae=a\}$ is defined by the equality
$$
\phi_e(a)=\phi(a)\,\,\,\forall a\in eMe, a\geq 0.
$$
One has the formula
\begin{equation}
\label{S(M)} S(M)=\bigcap_{e\not=0}Sp\Delta_{\phi_{e}},
\end{equation}
where $e$
varies over the nonzero projection of $M_\phi$ (see\cite{Con94}
p.472).
\begin{lem}
\label{M-phi} The von Neumann algebra  $M_\phi$ is trivial.
\end{lem}
In this case
$$S(M)=Sp\Delta =[0,\,\infty),$$
so the von Neumann algebra ${\mathfrak A}^{L,b}$ (and hence
algebra ${\mathfrak A}^{R,b}$) is type ${\rm III}_1$ factor.
\qed\end{pf}
\begin{pf*}{Proof of  Lemma \ref{M-phi}.}
We show that
\begin{equation}
\label{M=()'} M_\phi=(\Delta^{it},\,T^{R,b}_s\mid t\in {\mathbb
R},\,s \in B_0^{\mathbb N})^{\prime}.
\end{equation}
So $M_\phi$ is trivial means that the set of operators
\begin{equation}
\label{M'-irr} (\Delta^{it} ,\,T^{R,b}_s \mid t\in {\mathbb R},\,s
\in B_0^{\mathbb N})
\end{equation}
 is {\bf irreducible}. To prove (\ref{M=()'})
we get
$$
M_\phi=(a \in {\mathfrak A}^{L,b}\mid \Delta^{it}a=a\Delta^{it},
\,\,\forall t\in {\mathbb R })=(\Delta^{it}\mid t\in{\mathbb R
})'\cap{\mathfrak A}^{L,b}
$$
$$
=(\Delta^{it}\mid t\in{\mathbb R })'\cap({\mathfrak A}^{R,b})'
=(\Delta^{it}\mid t\in{\mathbb R })'\cap(T^{R,b}_s \mid s \in
B_0^{\mathbb N})^{\prime}=
$$
$$
(\Delta^{it},\,T^{R,b}_s\mid t\in {\mathbb R},\,s \in B_0^{\mathbb
N})^{\prime}
$$
{\bf Definition.} Recall (c.f. e.g. \cite{Dix69W}) that a non
necessarily bounded  self-adjoint operator $A$ in a Hilbert space
$H$ is said to be {\it affiliated} with a von Neumann algebra $M$
of operators in this Hilbert space $H$, if $\exp(itA)\in M$ for
all $t\in{\mathbb R}$. One then writes $A\,\,\eta\,\,M$.

 To prove
the irreducibility of $(\Delta^{it} ,\,T^{R,b}_s \mid t\in
{\mathbb R},\,s \in B_0^{\mathbb N})$ it is sufficient to prove
(see \cite{Kos92} p.258) that operators $f(x)\mapsto x_{kn}f(x)$
of multiplication in the space $H_b$ by the independent variables
$x_{kn}$ are affiliated to the von Neumann algebra
$$
(M_\phi)'=(\Delta^{it},\,T^{R,b}_s\mid t\in {\mathbb R},\,s \in
B_0^{\mathbb N})''.
$$
In this case the operator $A$ commuting with $\Delta^{it}$ and
$T^{R,b}_s$ is operator of multiplication by some function $a(x)$.
If we use commutation relation $[A,T^{R,b}_s]=0,\,\,s \in
B_0^{\mathbb N}$ we obtain $a(x)=a(xs)\,\,{\rm mod}\mu$. Using the
ergodocity of the measure $\mu_b$ with respect of the right action
of the group $B_0^{\mathbb N}$ we conclude that $a(x)={\rm
const}\,\,{\rm mod}\mu$ i.e. $A$ is scalar operator.

If we denote
$$A^R_{kn}=(d/dt)T^{R,b}_{I+tE_{kn}}\mid_{t=o}$$
we have (see for example \cite{Kos90,Kos92,Kos94})
\begin{equation}
\label{A^R-kn}
A^{R}_{kn}=\sum_{r=1}^{k-1}x_{kr}D_{rn}+D_{kn},\quad 1\leq k<n.
\end{equation}
The direct calculation shows that
\begin{equation}
\label{[13[23,D]]} [A^R_{13},[A^R_{23}
,\,\ln\Delta]]=2b_{13}x_{12},
\end{equation}
\begin{equation}
\label{[12[23,D]]} [A^R_{12},[A^R_{23}
,\,\ln\Delta]]=2b_{13}x_{13}.
\end{equation}

{\bf Idea: to obtain in a similar way all variables $x_{kn}.$}

 Let us denote by $X^{-1}$ the inverse matrix to the
upper triangular matrix $X=I+x=I+\sum_{k<n}x_{kn}E_{kn}\in
B^{\mathbb N}$
$$
X^{-1}=(I+x)^{-1}=I+\sum_{k<n}x_{kn}^{-1}E_{kn}\in B^{\mathbb N}.
$$
We have by definition $X^{-1}X=XX^{-1}=I$ hence
\begin{equation}
\label{x_{kn}{-1}} \left(XX^{-1}\right)_{kn}=
\sum_{r=k}^{n}x_{kr}x_{rn}^{-1}=\delta_{kn}=\sum_{r=k}^{n}x_{kr}^{-1}x_{rn}
=\left(X^{-1}X\right)_{kn} ,\quad k\leq n,
\end{equation}
hence
$$
x_{kn}^{-1}+\sum_{r=k+1}^{n-1}x_{kr}x_{rn}^{-1}+x_{kn}=0=
x_{kn}+\sum_{r=k+1}^{n-1}x_{kr}x_{rn}^{-1}+x_{kn}^{-1},\quad k< n,
$$
and
\begin{equation}
\label{x{kn}(-1)1}
x_{kn}^{-1}=-x_{kn}-\sum_{r=k+1}^{n-1}x_{kr}x_{rn}^{-1}=
-x_{kn}-\sum_{r=k+1}^{n-1}x_{kr}^{-1}x_{rn}.
\end{equation}
We can write also
\begin{equation}
\label{x{kn}(-1)2} x_{kn}^{-1}=-\sum_{r=k+1}^{n}x_{kr}x_{rn}^{-1}=
-\sum_{r=k}^{n-1}x_{kr}^{-1}x_{rn}.
\end{equation}
There is also the explicit formula for $x_{kn}^{-1}$ (see
\cite{Kos88} formula (4.4)) $x_{kk+1}^{-1}=-x_{kk+1}$
\begin{equation}
\label{x{kn}(-1)}
x_{kn}^{-1}=-x_{kn}+\sum_{r=1}^{n-k-1}(-1)^{r+1}\sum_{k\leq
i_1<i_2<...<i_r\leq n }x_{ki_1}x_{i_1i_2}...x_{i_rn},\quad k<n-1.
\end{equation}
\begin{rem}
 Using (\ref{x{kn}(-1)}) we see that $x_{kn}^{-1}$ depends only on $x_{rs}$ with $k\leq r<s\leq n$.
\end{rem}
 Using (\ref{x{kn}(-1)2}) we have
\begin{equation}
\label{(.+.)(.-.)}
x_{kn}+x_{kn}^{-1}=-\sum_{r=k+1}^{n-1}x_{kr}x_{rn}^{-1},\quad
x_{kn}-x_{kn}^{-1}=2x_{kn}-\sum_{r=k+1}^{n-1}x_{kr}x_{rn}^{-1}.
\end{equation}
Let us denote
\begin{equation}
\label{w_kn}
w_{kn}:=w_{kn}(x):=(x_{kn}+x_{kn}^{-1})(x_{kn}-x_{kn}^{-1}).
\end{equation}
 Using (\ref{mu-b}) we get
\begin{equation}
\label{Delta(x)} \Delta(x)=\frac{d\mu_b(x)}{d\mu_b(x^{-1})}
=\exp\left[-\sum_{k<n}b_{kn}\left(x_{kn}^2-(x_{kn}^{-1})^2\right)
\right]=\exp\left[-\sum_{k<n}b_{kn}w_{kn}(x)\right].
\end{equation}
$$
-\ln\Delta(x)=\sum_{k<n}b_{kn}\left[x_{kn}^2-(x_{kn}^{-1})^2\right]=
\sum_{k<n}b_{kn}(x_{kn}+x_{kn}^{-1})(x_{kn}-x_{kn}^{-1})
$$
$$
\sum_{k<n}b_{kn}(x_{kn}+x_{kn}^{-1})[2x_{kn}-(x_{kn}+x_{kn}^{-1})]=
\sum_{k<n}b_{kn}w_{kn}(x).
$$
To study the action of the operators
$A^{R}_{kn}=\sum_{r=1}^{k-1}x_{rk}D_{rn}+D_{kn}$ on the function
$\ln\Delta(x)$ we need to know the action of $D_{pq}$ on
$x_{kn}^{-1}$.
\begin{lem}
\label{} We have
\begin{equation}
\label{[D,x{-1}]} [D_{pq},x_{kn}^{-1}]= \left\{\begin{array}{ll}
-x_{kp}^{-1}x_{qn}^{-1},&\text{\,if\,}\,k\leq p<q\leq n,\\
0,&\text{\,otherwise\,}.
\end{array}\right.
\end{equation}
\end{lem}
\begin{pf} We prove (\ref{[D,x{-1}]}) by induction in $p:k\leq p<q\leq
n$. For  $p=k$ using (\ref{(.+.)(.-.)}) we have
$$
[D_{kq},x_{kn}^{-1}]=-[D_{kq},x_{kn}+\sum_{r=k+1}^{n-1}x_{kr}x_{rn}^{-1}]=
-[D_{kq},x_{kq}x_{qn}^{-1}]=-x_{qn}^{-1}=-x_{kk}^{-1}x_{qn}^{-1},
$$
so (\ref{[D,x{-1}]}) holds for $p=k$.

Let us suppose that (\ref{[D,x{-1}]}) holds for all $(p,q)$ with
$k\leq p<s\leq n,\,\,k\leq p<q\leq n$. We prove that than
(\ref{[D,x{-1}]}) holds also for $(s,q)\,:\,s<q\leq n$. Indeed we
have
$$
[D_{sq},x_{kn}^{-1}]=-[D_{sq},x_{kn}^{-1}+\sum_{r=k+1}^{n-1}x_{kr}x_{rn}^{-1}]=
-\sum_{r=k+1}^{s}x_{kr}[D_{sq},x_{rn}^{-1}]
$$
$$
=\sum_{r=k+1}^{s}x_{kr}x_{rs}^{-1}x_{qn}^{-1}
\stackrel{(\ref{x{kn}(-1)1} )}{=} x_{ks}^{-1}x_{qn}^{-1}.
$$
\qed\end{pf}

Using (\ref{[D,x{-1}]}) we get
\begin{equation}
\label{[D,x+x{-1}]} [D_{pq},x_{kn}+x_{kn}^{-1}]=
\left\{\begin{array}{ll}
-x_{kp}^{-1}x_{qn}^{-1},&\text{\,if\,}\,k\leq p<q\leq n,\,(p,q)\not=(k,n)\\
0,&\text{\,otherwise\,}.
\end{array}\right.
\end{equation}
Using (\ref{[D,x+x{-1}]}) we have
\begin{equation}
\label{[D,w]} [D_{pq},(x_{kn}+x_{kn}^{-1})(x_{kn}-x_{kn}^{-1})]=
\left\{\begin{array}{ll}
2\,x_{kp}^{-1}x_{qn}^{-1}x_{kn}^{-1},&\text{\,if\,}\,k\leq p<q\leq n,\,(p,q)\not=(k,n)\\
2(x_{kn}+x_{kn}^{-1}),&\text{\,if\quad}(p,q)=(k,n)\\
0,&\text{\,otherwise\,}.
\end{array}\right.
\end{equation}
Indeed, if $k\leq p<q\leq n,\,(p,q)\not=(k,n)$ we have
$$
[D_{pq},(x_{kn}+x_{kn}^{-1})(x_{kn}-x_{kn}^{-1})]=
[D_{pq},(x_{kn}+x_{kn}^{-1})(2x_{kn}-(x_{kn}+x_{kn}^{-1}))]
$$
$$
=[D_{pq},(x_{kn}+x_{kn}^{-1})](2x_{kn}-(x_{kn}+x_{kn}^{-1}))-
(x_{kn}+x_{kn}^{-1})[D_{pq},(x_{kn}+x_{kn}^{-1})]=
$$
$$
-2x_{kn}^{-1}[D_{pq},(x_{kn}+x_{kn}^{-1})]
\stackrel{(\ref{[D,x+x{-1}]} )}{=}
2x_{kp}^{-1}x_{qn}^{-1}x_{kn}^{-1}.
$$
\begin{lem}
\label{l.[A^R,w]}  We have
\begin{equation}
\label{[A^R,w]}
[A^R_{mm+1} ,\,w_{kn}]= \left\{\begin{array}{ll}
0,                        &\text{\,if\,}\,\,\,k<n\leq m\\
2x_{km}x_{km+1}           &\text{\,if\,}\,\,\,n=m+1,\,\,1\leq k\leq m-1\\
0,                        &\text{\,if\,}\,\,\,1\leq k\leq m-1,\,\,m+1<n\\
2x_{mn}^{-1}x_{m+1n}^{-1},&\text{\,if\,}\,\,\,k=m,\,\, n\geq m+2\\
0,                        &\text{\,if\,}\,\,\,\,\, m+1\leq k<n.
\end{array}\right.
\end{equation}
hence
\begin{equation}
\label{[A^R,lnDe]} - [A^R_{mm+1}
,\,\ln\Delta]=2\sum_{r=1}^{m-1}b_{rm+1}x_{rm}x_{rm+1} +
2\sum_{n=m+2}^\infty b_{mn}x_{mn}^{-1}x_{m+1n}^{-1}.
\end{equation}
\end{lem}
\begin{pf}
Since $$ A^R_{mm+1}=\sum_{r=1}^{m-1}x_{rm}D_{rm+1}+D_{mm+1}
$$
and $w_{kn},\,\,k<n\leq m$ do not depend on $x_{rm+1},\,\,1\leq
r\leq m+1$  we conclude that $[A^R_{mm+1} ,\,w_{kn}]=0$ for
$k<n\leq m$ and  $m+1\leq k<n$.

Let $n=m+1$, since $[D_{rm+1},\,w_{km+1}]=0$ for $1\leq r<k$ we
get
$$
[A^R_{mm+1} ,\,w_{km+1}] =
\sum_{r=k}^{m-1}x_{rm}[D_{rm+1},\,w_{km+1}]+[D_{mm+1},\,w_{km+1}]=
$$
$$
2\left(x_{km}(x_{km+1}+x_{km+1}^{-1})+
\sum_{r=k+1}^{m-1}x_{rm}x_{kr}^{-1}x_{km+1}^{-1}+x_{km}^{-1}x_{km+1}^{-1}\right)=
$$
$$
2\left(x_{km}x_{km+1}+\left(x_{km}+\sum_{r=k+1}^{m-1}x_{kr}^{-1}x_{rm}+x_{km}^{-1}
\right)x_{km+1}^{-1}\right) \stackrel{(\ref{x{kn}(-1)1})}{=}
2x_{km}x_{km+1}.
$$
Similarly, for $1\leq k\leq m-1,\,\,m+1<n$ we get
$$
[A^R_{mm+1} ,\,w_{kn}] =
\sum_{r=k}^{m-1}x_{rm}[D_{rm+1},\,w_{kn}]+[D_{mm+1},\,w_{kn}]=
$$
$$
2\left(x_{km}x_{m+1n}^{-1}+
\sum_{r=k+1}^{m-1}x_{rm}x_{kr}^{-1}x_{m+1n}^{-1}+x_{km}^{-1}x_{m+1n}^{-1}\right)
$$
$$
2\left(x_{km}+
\sum_{r=k+1}^{m-1}x_{rm}x_{kr}^{-1}+x_{km}^{-1}\right)x_{m+1n}^{-1}\stackrel{(\ref{x{kn}(-1)1})}{=}0.
$$
Finally if $k=m$ and $n\geq m+2$ we have as before
$$
[A^R_{mm+1} ,\,w_{mn}]=[D_{mm+1},\,w_{mn}]
\stackrel{(\ref{[D,w]})}{=}2x_{mn}^{-1}x_{m+1n}^{-1}.
$$
\qed\end{pf}

We consider the action of $A^R_{mm+1}$ on $\ln\Delta$.

Let $m=2$. Since
$$[A^R_{23} ,\,w_{13}]=2b_{13}x_{12}x_{13},\,\,\,
[A^R_{23} ,\,w_{1n}]=0,\,\, n\geq 4,\,\,\,[A^R_{23}
,\,w_{kn}]=0,\,\,3\leq k<n,
$$
 we have
$$
-[A^R_{23} ,\,\ln\Delta]=2b_{13}x_{12}x_{13}+2\sum_{n=4}^\infty
b_{2n}x_{2n}^{-1}x_{3n}^{-1},
$$
hence
$$
-[A^R_{12},[A^R_{23} ,\,\ln\Delta]]=2b_{13}x_{13},
$$
$$
-[A^R_{13},[A^R_{23} ,\,\ln\Delta]]=2b_{13}x_{12}.
$$
The last two equations gives us $x_{12},x_{13}\,\eta\, {\mathfrak
A}$.

Let $m=3$. Since
$$
[A^R_{34} ,\,w_{13}]=0,\,\,[A^R_{34} ,\,w_{14}]=2x_{13}x_{14},\,\,
[A^R_{34} ,\,w_{24}]=2x_{23}x_{24},
$$
$$
[A^R_{34} ,\,w_{1n}]=[A^R_{34} ,\,w_{1n}]=0,\,\, [A^R_{34}
,\,w_{3n}]=b_{3n}x_{3n}^{-1}x_{4n}^{-1}, n\geq 5,
$$
$$
[A^R_{34} ,\,w_{kn}]=0,\,\,4\leq k<n,
$$
we have
$$
-[A^R_{34} ,\,\ln\Delta]=2b_{14}x_{13}x_{14}+2b_{24}x_{23}x_{24}+
2\sum_{n=5}^\infty b_{3n}x_{3n}^{-1}x_{4n}^{-1},
$$
hence
$$
-[A^R_{23},[A^R_{34}
,\,\ln\Delta]]=2b_{14}x_{12}x_{14}+2b_{24}x_{24}
$$
$$
-[A^R_{12}[A^R_{23},[A^R_{34} ,\,\ln\Delta]]]=2b_{14}x_{14},
$$
$$
-[A^R_{24},[A^R_{34}
,\,\ln\Delta]]=2[x_{12}D_{14}+D_{24},b_{14}x_{13}x_{14}+b_{24}x_{23}x_{24}]=
2b_{14}x_{12}x_{13}+2b_{24}x_{23},
$$
Since $x_{12},x_{13}\,\eta\, {\mathfrak A}$ from the latter
equation we conclude that $x_{23}\,\eta \,{\mathfrak A}$. The
previous equation gives us $x_{14}\,\eta \,{\mathfrak A}$ and the
equation before gives $x_{24}\,\eta \,{\mathfrak A}$. Finally we
conclude that $x_{14},x_{24},x_{23}\,\eta \,{\mathfrak A}$.

Let us suppose that we have obtained the variables $x_{rm},1\leq
r\leq m-2$ and $x_{m-2,m-1}$. We prove that we can obtain the
following variables $x_{rm+1},1\leq r\leq m-1$ and $x_{m-1m}.$

Indeed we calculate the action of the following sequence of
operators on the result: $A^R_{m-1,m},\,A^R_{m-2,m-1}$ etc. till
$A^R_{12}$. We obtain
$$
-[A^R_{m-1,m},[A^R_{mm+1} ,\,\ln\Delta]]=
2\left(\sum_{r=1}^{m-2}b_{r,m+1}x_{r-1,m}x_{r,m+1}+b_{m-1,m+1}x_{m-1,m+1}\right),
$$
$$
-[A^R_{m-2,m-1},[A^R_{m-1,m},[A^R_{mm+1} ,\,\ln\Delta]]]
$$
$$
=
2\left(\sum_{r=1}^{m-3}b_{r,m+1}x_{r-2,m}x_{r,m+1}+b_{m-2,m+1}x_{m-2,m+1}\right),
$$
$$
-[A^R_{m-s,m-s+1},[A^R_{m-s+1,m-s+2},...[A^R_{m-1,m}[A^R_{mm+1}
,\,\ln\Delta]]...]]
$$
$$
=
2\left(\sum_{r=1}^{m-s-1}b_{r,m+1}x_{r,m-s}x_{r,m+1}+b_{m-s,m+1}x_{m-s,m+1}\right),\,\,1\leq
s\leq m,
$$
$$
-[A^R_{34},...[A^R_{mm+1}
,\,\ln\Delta]...]=2(b_{1,m+1}x_{13}x_{1,m+1}+b_{2,m+1}x_{23}x_{2,m+1}+b_{3,m+1}x_{3,m+1}),
$$
$$
-[A^R_{23},[A^R_{34},...[A^R_{mm+1}
,\,\ln\Delta]...]]=2(b_{1,m+1}x_{12}x_{1,m+1}+b_{2,m+1}x_{2,m+1}),
$$
$$
-[A^R_{12},[A^R_{23},[A^R_{34},...[A^R_{mm+1}
,\,\ln\Delta]...]]]=2b_{1,m+1}x_{1,m+1}.
$$
From the latter equation we conclude that $x_{1,m+1}\,\,\eta
\,\,{\mathfrak A}$. The last but one equation gives us
$x_{2,m+1}\,\,\eta\, \,{\mathfrak A}$ (since
$x_{12},x_{1,m+1}\,\,\eta\, \,{\mathfrak A}$) etc. i.e :
$x_{rm+1}\,\,\eta \,\,{\mathfrak A},\,\,1\leq r \leq m-1$.
$$
-[A^R_{m-1m+1}, [A^R_{mm+1}
,\,\ln\Delta]]=[\sum_{r=1}^{m-2}x_{rm-1}D_{rm+1}+D_{m-1m+1},
2\sum_{r=1}^{m-1}b_{rm+1}x_{rm}x_{rm+1}]=
$$
$$
2\sum_{r=1}^{m-2}b_{rm+1}x_{rm-1}x_{rm}+b_{m-1,m+1}x_{m-1,m},
$$
since $x_{rm-1},x_{rm}\,\,\eta \,\,{\mathfrak A}$ for $1\leq r\leq
m-2$ hence $x_{m-1,m}\,\,\eta \,\,{\mathfrak A}$. \qed\end{pf*}
To be sure  that all this argument works we should prove that all
involved operators are affiliated to the von Neumann algebra
$M_\phi'$ defined by (7). For example if  $A_{23}^R$ and $\Delta$
(and hence $\ln\Delta$) are affiliated to the von Neumann algebra
$M_\phi'$, why the operator $[A_{23}^R,\ln\Delta]$ is also
affiliated. In general, why the operators
$[A^R_{12},[A^R_{23},[A^R_{34},...[A^R_{mm+1} ,\,\ln\Delta]...]]]$
are affiliated?
\begin{rem}
\label{[A,B]in M} In general we do not know whether the commutator
$[A,B]$ of two operators $A$ and $B$ affiliated to the von Neumann
algebra is also affiliated.
\end{rem}
This is the reason, why we use another approach to prove that the
algebra $M_\phi$ is trivial.
\section{The von Neumann algebra  $M_\phi$ is trivial}
Since $M_\phi=(\Delta^{it},\,T^{R,b}_s\mid t\in {\mathbb R},\,s
\in B_0^{\mathbb N})^{\prime}$ (see (\ref{M=()'})) it is
sufficient to prove that the set of operators
$$
(\Delta^{is} ,\,T^{R,b}_t \mid s\in {\mathbb R},\,t \in
B_0^{\mathbb N})\subset M_\phi'
$$
is irreducible.

{\bf Idea of the proof}. {\it We show that the von Neumann
subalgebra in the algebra $M_\phi'$, generated by the following
operators
\begin{equation}
\label{max-ab}
(\{T^R_{t_{n}},\{T^R_{t_{n-1}},...\{T^R_{t_{1}},\Delta^{is}\}...\}\}\mid
s\in{\mathbb R},\,\,t_1,...,t_n\in B_0^{\mathbb N}),
\end{equation}
where $\{a,b\}:=aba^{-1}b^{-1}$ is {\bf the maximal abelian
subalgebra}. More precisely we prove that this subalgebra contains
all functions $\exp(isx_{kn}),\,\,k<n,\,\,s\in{\mathbb R}$}.

{\bf To prove the irreducibility of the algebra $M_\phi'$} (see
proof of the Lemma 14) we observe that if an bounded operator
commute with all $\exp(isx_{kn}),\,\,k<n,\,\,s\in{\mathbb R}$ then
this operator itself is an operator of multiplication by some
essentially bounded function $A=a(x)$. Commutation relation
$[T^{R,b}_t,A]=0$ for all $t\in B_0^{\mathbb N}$ gives us
$a(xt)=a(x)\,{\rm mod}\mu_b$ for all $t$. Since the measure
$\mu_b$ is {\bf $B_0^{\mathbb N}-$right ergodic} we conclude that
$A$ is trivial i.e. $A=a(x)=CI$.

We note that expressions in (\ref{max-ab}) are the "right" analog
of the left hand side of the expressions (\ref{[13[23,D]]}) and
(\ref{[12[23,D]]})
$$
[A^R_{13},[A^R_{23} ,\,\ln\Delta]]=2b_{13}x_{12},
$$
$$
[A^R_{12},[A^R_{23} ,\,\ln\Delta]]=2b_{13}x_{13},
$$
involving generators $A^R_{kn}$. In general, if we have two
subgroups of unitary operators $U(t)$ and $V(s)$ with the
generators $A$ and $B$, to obtain the commutator $[iA,iB]$ it is
sufficient to differentiate the following expression
$U(t)V(s)U(-t)$:
$$
\frac{\partial}{\partial t}\frac{\partial}{\partial
s}U(t)V(s)U(-t)\mid_{t=s=0}=[iA,iB].
$$
Indeed we have
$$
\frac{\partial}{\partial
s}U(t)V(s)U(-t)=U(t)iBV(s)U(-t),\quad\frac{\partial}{\partial
t}U(t)iBV(s)U(-t)\mid_{t=s=0}=
$$
$$
(iAU(t)iBV(s)U(-t)-U(t)iBV(s)iAU(-t))\mid_{t=s=0}=[iA,iB].
$$
We show that more convenient analog of the commutator $[iA,iB]$ is
{\it commutator} (in the group sence) of two one-parameter groups
$$\{U(t),V(s)\}:=U(t)V(s)U(t)^{-1}V(s)^{-1}=U(t)V(s)U(-t)V(-s).$$
\begin{lem}
\label{l.[T_t,g]} For the operator $g$ of multiplication on the
function $g:f(x)\mapsto g(x)f(x)$ in the space $H_b=L_2(B^{\mathbb
N},d\mu_b)$ we have
$$
T^R_tg(x)T^R_{t^{-1}}=g(xt),\,\,t\in B_0^{\mathbb N}.
$$
\end{lem}
\begin{pf} We have
$$
f(x)\stackrel{g(x)T^R_{t^{-1}}}{\mapsto}
g(x)\left(\frac{d\mu(xt^{-1})}{d\mu(x)}\right)^{1/2}f(xt^{-1})
\stackrel{T^R_{t}}{\mapsto}
$$
$$
\left(\frac{d\mu(xt)}{d\mu(x)}\right)^{1/2}g(xt)\left(\frac{d\mu(x)}
{d\mu(xt)}\right)^{1/2}f(x)= g(xt)f(x).
$$
\qed\end{pf} Using the lemma we have
$$
T^R_t\Delta^{is}(x)T^R_{t^{-1}}=\Delta^{is}(xt).
$$
Using (\ref{Delta(x)}) we have
$$
\Delta^{is}(x)=\exp\left( -is\sum_{k+1<n}
b_{kn}(x_{kn}+x_{kn}^{-1})[2x_{kn}-(x_{kn}+x_{kn}^{-1})]\right)=
$$
\begin{equation}
\label{Del^{is}(x)} \exp\left( -is\sum_{k+1<n} b_{kn}w_{kn}(x)
\right),
\end{equation}
where
$w_{kn}(x)\!=\!(x_{kn}+x_{kn}^{-1})[2x_{kn}-(x_{kn}+x_{kn}^{-1})]$
(see (\ref{w_kn})).

We would like {\bf to obtain the functions} $\exp(isx_{kn})$ using
the expressions (\ref{max-ab}). To simplify the situation we
consider firstly the projections of all considered object: the
measure $\mu_b^{(k)}$, the generators $A^{R,(k)}_{kn}$, operator
$\Delta_{(k)}$ algebra $M^{(k)}:=(M_\phi')^{(k)}$ etc. on the
following subspace $X^{(k)},\,\,k\geq 2$ of the space $B^{\mathbb
N}$:
$$ X^{(2)}=\left(\begin{smallmatrix}
1&x_{12}&x_{13}&...&x_{1n}&...\\
0&1     &x_{23}&...&x_{2n}&...\\
\end{smallmatrix}\right),\quad
X^{(3)}=\left(\begin{smallmatrix}
1&x_{12}&x_{13}&x_{14}&...&x_{1n}&...\\
0&1     &x_{23}&x_{14}&...&x_{2n}&...\\
0&0     &1     &x_{34}&...&x_{3n}&...\\
\end{smallmatrix}\right),\,\,\text{\rm etc}.
$$
Note that
\begin{equation}
\label{(X^2)^{-1}}
 \left(\begin{smallmatrix}
1&x_{12}&x_{13}&...&x_{1n}&...\\
0&1     &x_{23}&...&x_{2n}&...\\
\end{smallmatrix}\right)^{-1}=
\left(\begin{smallmatrix}
1&-x_{12}&-x_{13}+x_{12}x_{23}&...&-x_{1n}+x_{12}x_{2n}&...\\
0&1     &-x_{23}&...&-x_{2n}&...\\
\end{smallmatrix}\right).
\end{equation}
We have for the corresponding projections on $X^{(2)}$:
$$
A^{R}_{1n}=D_{1n},\quad A^{R}_{2n}=x_{12}D_{1n}+D_{2n},\quad
A^{R,(2)}_{kn}=x_{1k}D_{1n}+x_{2k}D_{2n},\,\,2<k<n,
$$
$$
w_{1n}(x)\!=\!(x_{1n}+x_{1n}^{-1})(x_{1n}-x_{1n}^{-1})=
x_{12}x_{2n}(2x_{1n}-x_{12}x_{2n}),\,\,w_{2n}(x)=0,
$$
hence
$$
\Delta_{(2)}^{is}(x):=\exp\left( -is\sum_{k=3}^\infty
b_{1n}w_{1n}(x) \right)=\exp\left( -is\sum_{k=3}^\infty
b_{1n}x_{12}x_{2n}(2x_{1n}-x_{12}x_{2n}) \right).
$$
Let us denote by
\begin{equation}
\label{E_kn(t)} E_{kn}(t):=I+tE_{kn},\quad T_{kn}(t)=
T^R_{E_{kn}(t)},\,\, k<n,\,t\in{\mathbb R}
\end{equation} the corresponding one-parameter subgroups.
We have
$$
\left(\begin{smallmatrix} x_{12}&x_{1m}\\
1&x_{2m}\\
\end{smallmatrix}\right)
\stackrel{E_{2m}(t) }{\mapsto}
\left(\begin{smallmatrix} x_{12}&x_{1m}+tx_{12}\\
1&x_{2m}+t\\
\end{smallmatrix}\right),\,\,w_{1n}(xE_{2m}(t))=
\left\{\begin{smallmatrix} w_{1n}(x)&\text{\,\,if\,\,}n\not=m\\
w_{1m}(xE_{2m}(t))&\text{\,\,if\,\,}n=m\\
\end{smallmatrix}\right.
$$
so  using Lemma \ref{l.[T_t,g]} we get
$$
\{T_{2m}(t),\Delta_{(2)}^{is}(x)\}=T_{2m}(t)\Delta_{(2)}^{is}(x)
T_{2m}(-t)\Delta_{(2)}^{-is}(x)=\Delta_{(2)}^{is}(xE_{2m}(t))\Delta_{(2)}^{-is}(x)=
$$
$$
\exp\left( -is\left[\sum_{k=3,k\not=m}^\infty
b_{1n}w_{1n}(x)+b_{1m}w_{1m}(xE_{2m}(t))\right] \right)\exp\left(
is\sum_{k=3}^\infty b_{1n}w_{1n}(x) \right)=
$$
$$
\exp\left( -isb_{1m}[w_{1m}(xE_{2m}(t))-w_{1m}(x)] \right)=
\exp\left( isb_{1m}(2tx_{12}x_{1m}+t^2x_{12}^2)\right),
$$
since
$$
w_{1m}(xE_{2m}(t))-w_{1m}(x)=x_{12}(x_{2m}+t)[2(x_{1m}+tx_{12})-x_{12}(x_{2m}+t)]-
$$
$$
x_{12}x_{2m}(2x_{1m}-x_{12}x_{2m})=x_{12}[tx_{12}x_{2m}+t(2x_{1m}-x_{12}x_{2m})+t^2x_{12}]
=2tx_{12}x_{1m}+t^2x_{12}^2.
$$
Let us denote
\begin{equation}
\label{[2m,Del]}
\phi_{t,s}(x):=\{T_{2m}(t),\Delta_{(2)}^{is}(x)\}=\exp\left(
isb_{1m}(2tx_{12}x_{1m}+t^2x_{12}^2)\right).
\end{equation}
Using Lemma \ref{l.[T_t,g]} we get
$$
\{T_{1m}(t_1),\{T_{2m}(t),\Delta_{(2)}^{is}(x)\}\}=\{T_{1m}(t_1),\phi_{t,s}(x)\}=
$$
$$
T_{1m}(t_1)\phi_{t,s}(x)T_{1m}(-t_1)(\phi_{t,s}(x))^{-1}=
\phi_{t,s}(xE_{1m}(t_1))(\phi_{t,s}(x))^{-1}=
$$
$$
\exp\left[ isb_{1m}(2tx_{12}(x_{1m}+t_1)+t^2x_{12}^2)-
isb_{1m}(2tx_{12}x_{1m}+t^2x_{12}^2)\right]=
$$
$$
\exp\left( isb_{1m}x_{12}2tt_1\right).
$$
Finally we get for $X^{(2)}$
$$
\exp(isx_{12})\in M^{(2)}:=(M_\phi')^{(2)}.
$$
Using (\ref{[2m,Del]}) we conclude that
$$
\exp(isx_{12}x_{1m})\in M^{(2)}.
$$
Applying again $T_{12}(t)$ and $T_{1m}(t)$ we get
$$
\{T_{12}(t),\exp(isx_{12}x_{1m})\}= T_{12}(t)\exp(isx_{12}x_{1m})
T_{12}(-t)
\exp(-isx_{12}x_{1m})=
$$
$$
\exp(is(x_{12}+t)x_{1m}-isx_{12}x_{1m})= \exp(istx_{12}),
$$
$$
\{T_{1m}(t),\exp(isx_{12}x_{1m})\}=T_{1m}(t)\exp(isx_{12}x_{1m})
T_{1m}(-t),\exp(-isx_{12}x_{1m})=
$$
$$
\exp(isx_{12}(x_{1m}+t)-isx_{12}x_{1m})= \exp(istx_{1m}).
$$
At last we conclude that for $X^{(2)}$ we have
$\exp(isx_{12}),\,\,\exp(isx_{1m})\in M^{(2)}$ in particular
\begin{equation}
\label{x_{12},x_{13},(2)}
 \exp(isx_{12}),\,\,\exp(isx_{13})\in
M^{(2)}.
\end{equation}
For $X^{(3)}$ and the corresponding projections we have
$$
\left(\begin{smallmatrix}
1&x_{12}&x_{13}&x_{14}&...&x_{1n}&...\\
0&1     &x_{23}&x_{14}&...&x_{2n}&...\\
0&0     &1     &x_{34}&...&x_{3n}&...\\
\end{smallmatrix}\right)^{-1}=
$$
$$
\left(\begin{smallmatrix}
1&-x_{12}&-x_{13}+x_{12}x_{23}&-x_{14}+x_{12}x_{24}+x_{13}x_{34}+x_{12}x_{23}x_{34}&
...&-x_{1n}+x_{12}x_{2n}+x_{13}x_{3n}+x_{12}x_{23}x_{3n}&...\\
0&1     &-x_{23}&-x_{24}+x_{23}x_{34}&...&-x_{2n}+x_{23}x_{3n}&...\\
0&0     &1     &-x_{34}&...&x_{3n}&...\\
\end{smallmatrix}\right)=
$$
\begin{equation}
\label{(X^3)^{-1}} \left(\begin{smallmatrix}
1&-x_{12}&-x_{13}-x_{12}^{-1}x_{23}&-x_{14}-x_{12}^{-1}x_{24}-x_{13}^{-1}x_{34}&...&
-x_{1n}-x_{12}^{-1}x_{2n}-x_{13}^{-1}x_{3n}&...\\
0&1     &-x_{23}&-x_{24}-x_{23}^{-1}x_{34}&...&-x_{2n}-x_{23}^{-1}x_{3n}&...\\
0&0     &1     &-x_{34}&...&-x_{3n}&...\\
\end{smallmatrix}\right),
\end{equation}
$$
A^{R}_{1n}=D_{1n},\quad A^{R}_{2n}=x_{12}D_{1n}+D_{2n},\quad
A^{R}_{3n}=x_{13}D_{1n}+x_{23}D_{2n}+D_{3n},\,\,3<n.
$$
We have
$$\Delta^{is}_{(3)}(x)=\exp\left( -is\left[\sum_{n=3}^\infty
b_{1n}w_{1n}(x)+\sum_{n=4}^\infty b_{2n}w_{2n}(x) \right]\right)=
$$
$$
\exp\left( -is\left[\sum_{n=3}^\infty
b_{1n}(x_{1n}+x_{1n}^{-1})[2x_{1n}-(x_{1n}+x_{1n}^{-1})]\right]\right)\times
$$
$$
\exp\left( -is\left[\sum_{n=4}^\infty
b_{2n}(x_{2n}+x_{2n}^{-1})[2x_{2n}-(x_{2n}+x_{2n}^{-1})
\right]\right).
$$
By the same procedure as in the case of the space $X^{(2)}$ we can
obtain that
\begin{equation}
\label{x_{12},x_{13},(3)}
\exp(isx_{12}),\,\,\exp(isx_{13})\in
M^{(3)}.
\end{equation}
We show that
\begin{equation}
\label{[34,Del_3]} \{T_{34}(t),\Delta_{(3)}^{is}(x)\}= \exp\left(
is\left[b_{14}(2tx_{13}x_{14}+t^{2}x_{13}^2)+
b_{24}(2tx_{23}x_{24}+t^2x_{23}^2)\right]\right).
\end{equation}
(compare with (\ref{[2m,Del]})). Indeed we have
$$
\{T_{34}(t),\Delta_{(3)}^{is}(x)\}=
T_{34}(t)\Delta_{(3)}^{is}(x)T_{34}(-t)\Delta_{(3)}^{-is}(x)=
$$
$$
\Delta_{(3)}^{is}(xE_{34}(t))\Delta_{(3)}^{-is}(x)=
$$
$$
\exp\left( -is\left( b_{14}[w_{14}(xE_{34}(t))-w_{14}(x)]+
b_{24}[w_{24}(xE_{34}(t))-w_{24}(x)] \right)\right),
$$
which implies (\ref{[34,Del_3]}), since
$$
w_{14}(x)=(x_{14}+x_{14}^{-1})[2x_{14}-(x_{14}+x_{14}^{-1})]=
-(x_{12}^{-1}x_{24}+x_{13}^{-1}x_{34})[2x_{14}+x_{12}^{-1}x_{24}+x_{13}^{-1}x_{34}],
$$
and
$$
w_{14}(xE_{34}(t))-w_{14}(x)=
$$
$$
-[x_{12}^{-1}(x_{24}+tx_{23})+x_{13}^{-1}(x_{34}+t)]
[2(x_{14}+tx_{13})+x_{12}^{-1}(x_{24}+tx_{23})+x_{13}^{-1}(x_{34}+t)]
$$
$$
+(x_{12}^{-1}x_{24}+x_{13}^{-1}x_{34})[2x_{14}+x_{12}^{-1}x_{24}+x_{13}^{-1}x_{34}]=
$$
$$
-t\left[ (x_{12}^{-1}x_{24}+x_{13}^{-1}x_{34})
(2x_{13}+x_{12}^{-1}x_{23}+x_{13}^{-1}) +
(x_{12}^{-1}x_{23}+x_{13}^{-1})(2x_{14}+x_{12}^{-1}x_{24}+x_{13}^{-1}x_{34})
\right]
$$
$$
-t^2(x_{12}^{-1}x_{23}+x_{13}^{-1})(2x_{13}+x_{12}^{-1}x_{23}+x_{13}^{-1})=
-t[-(x_{14}+x_{14}^{-1})x_{13}-x_{13}(x_{14}+x_{14}^{-1})]+t^2x_{13}x_{13}=
$$
$$
2tx_{13}x_{14}+t^{2}x_{13}^2.
$$
Using (\ref{x_{12},x_{13},(3)}) and (\ref{[34,Del_3]}) we get
$$
\phi_{t,s}^{(3)}(x):=\exp\left( is\left[b_{14}2tx_{13}x_{14}+
b_{24}(2tx_{23}x_{24}+t^2x_{23}^2)\right]\right)\in M^{(3)},
$$
hence
$$
\{T_{13}(t_1),\phi_{t,s}^{(3)}(x)\}=T_{13}(t_1)\phi_{t,s}^{(3)}(x)
T_{13}(-t_1)(\phi_{t,s}^{(3)}(x))^{-1}=\exp\left(
istt_1b_{14}2tx_{14}\right),
$$
so $\exp(isx_{14})\in M^{(3)}$ and $\exp[ is
b_{24}(2tx_{23}x_{24}+t^2x_{23}^2)]\in M^{(3)}$. Similarly we get
$$
\{T_{24}(t_1),\exp[ is b_{24}(2tx_{23}x_{24}+t^2x_{23}^2)]\}
=\exp( is b_{24}tt_1x_{23}),
$$
so $\exp (isx_{23}),\exp(is x_{23}x_{24})\in M^{(3)}$. At last we
get
$$
\{T_{24}(t_1),\exp(is x_{23}x_{24})\}=\exp(ist_1x_{24}).
$$
Finally we can obtain $\exp(isx_{kn})$ in the following order on
the {\bf first step}:
$$
 \exp (isx_{12}),\,\,\exp(is x_{13});\,\,
$$
on the {\bf second step}:
$$
\exp(is x_{14}),\,\,\exp (isx_{23}),\,\,\exp(isx_{24})\in M^{(3)},
$$
or symbolically  in the following {\bf order}:
$$
\left(\begin{smallmatrix}
1&x_{12}&x_{13}&x_{14}\\
0&1     &x_{23}&x_{24}\\
0&0     &1     &\\
\end{smallmatrix}\right),\quad \left(\begin{smallmatrix}
0&1_1&2_1&1_2\\
0&0  &2_2&3_2\\
0&0  &0  &\\
\end{smallmatrix}\right).
$$
In general we get {\bf the order}
\begin{equation}
\label{order}
 \left(\begin{smallmatrix}
1&x_{12}&x_{13}&x_{14}&x_{15}&x_{16}&x_{17}\\
0&1     &x_{23}&x_{24}&x_{25}&x_{26}&x_{27}\\
0&0     &1     &x_{34}&x_{35}&x_{36}&x_{37}\\
0&0&0&1&x_{45}&x_{46}&x_{47}\\
0&0&0&0&1&x_{56}&x_{57}\\
0&0&0&0&0&1&\\
\end{smallmatrix}\right),\quad
\left(\begin{smallmatrix}
0&1_1 &2_1&1_2&1_3&1_4&1_5\\
0&0 &2_2&3_2&2_3&2_4&2_5\\
0&0 &0&3_3&4_3&3_4&3_5\\
0&0 &0&0&4_4&5_4&4_5\\
0&0 &0&0&0&5_5&\\
\end{smallmatrix}\right).
\end{equation}

This {\bf order is right  in the general case} (without any
projections on $X^{(k)}$). To obtain $\exp(isx_{12})$ and
$\exp(isx_{13})$ on the {\bf first step} we get by Lemma
\ref{l.[T_t,g]}
$$
\{T_{23}(t),\Delta^{is}(x)\}=
T_{23}(t)\Delta^{is}(x)T_{23}(-t)\Delta^{-is}(x)=\Delta^{is}(xE_{23}(t))\Delta^{-is}(x)=
$$
\begin{equation}
\label{[T_23,De]} \exp\left\{ -is\left(\sum_{n=3}^\infty
b_{1n}[w_{1n}(xE_{23}(t))-w_{1n}(x)]+ \sum_{n=4}^\infty
b_{2n}[w_{2n}(xE_{23}(t))-w_{2n}(x)]\right) \right\}.
\end{equation}
Now we shall calculate $w_{1n}(xE_{23}(t))-w_{1n}(x)$ and
$w_{2n}(xE_{23}(t))-w_{2n}(x)$. We have by (\ref{(.+.)(.-.)})
$$
x_{1n}+x_{1n}^{-1}=-\sum_{r=2}^{n-1}x_{1r}x_{rn}^{-1},\quad
x_{2n}+x_{2n}^{-1}=-\sum_{r=3}^{n-1}x_{2r}x_{rn}^{-1}
$$
so we conclude that for $n>3$ holds
$$
(x_{1n}+x_{1n}^{-1})^{E_{23}(t)}=-\left(\sum_{r=2}^{n-1}x_{1r}x_{rn}^{-1}
\right)^{E_{23}(t)}= -\left(x_{12}x_{2n}^{-1}+ x_{13}x_{3n}^{-1}+
\sum_{r=4}^{n-1}x_{1r}x_{rn}^{-1}\right)^{E_{23}(t)}=
$$
$$
-\left(x_{12}(-x_{2n}-[x_{23}+t]x_{3n}^{-1}-\sum_{r=4}^{n-1}x_{2r}x_{rn}^{-1})
+[x_{13}+tx_{12}]x_{3n}^{-1}+\sum_{r=4}^{n-1}x_{1r}x_{rn}^{-1}\right)=
$$
$$
-\left(\sum_{r=2}^{n-1}x_{1r}x_{rn}^{-1}-tx_{12}x_{3n}^{-1}+tx_{12}x_{3n}^{-1}\right)=
x_{1n}+x_{1n}^{-1}.
$$
For $n=3$ we get
$x_{13}+x_{13}^{-1}=-x_{12}x_{23}^{-1}=x_{12}x_{23}$ hence
$$
(x_{13}+x_{13}^{-1})^{E_{23}(t)}=(x_{12}x_{23})^{E_{23}(t)}=
$$
$$
x_{12}[x_{23}+t]=x_{12}x_{23}+tx_{12}=x_{13}+x_{13}^{-1}-tx_{12}^{-1}.
$$
Finally we conclude that
\begin{equation}
\label{} (x_{1n}+x_{1n}^{-1})^{E_{23}(t)}=
\left\{\begin{array}{ll}
x_{1n}+x_{1n}^{-1},&\text{\,if\,}\,3<n,\\
x_{13}+x_{13}^{-1}+tx_{12},&\text{\,if\,}\,\,n=3
\end{array}
\right.
\end{equation}
and
\begin{equation}
\label{} (x_{1n}\pm x_{1n}^{-1})^{E_{23}(t)}=
\left\{\begin{array}{ll}
x_{1n}\pm x_{1n}^{-1},&\text{\,if\,}\,3<n,\\
x_{13} \pm x_{13}^{-1}+ tx_{12},&\text{\,if\,}\,\,n=3
\end{array}
\right.
\end{equation}
since
$$
(x_{13}-x_{13}^{-1})^{E_{23}(t)}=(2x_{13}-(x_{13}+x_{13}^{-1}))^{E_{23}(t)}=
2[x_{13}+tx_{12}]-(x_{13} + x_{13}^{-1}+ tx_{12})
$$
$$
=x_{13}-x_{13}^{-1}+tx_{12}.
$$
We have $w_{1n}(xE_{23}(t))-w_{1n}(x)=0$ for $n>3$. For $n=3$
holds
$$
w_{13}(xE_{23}(t))-w_{13}(x)=(x_{13}+x_{13}^{-1}+tx_{12})
(x_{13}-x_{13}^{-1}+tx_{12})-(x_{13}+x_{13}^{-1})(x_{13}-x_{13}^{-1})
$$
$$
=tx_{12}(x_{13}+x_{13}^{-1}+x_{13}-x_{13}^{-1})+t^2x_{12}^2=
2tx_{12}x_{13}+t^2x_{12}^2.
$$
Finally
%
%
\begin{equation}
\label{w_1n^23} w_{1n}(xE_{23}(t))-w_{1n}(x)=
\left\{\begin{array}{ll}
0,&\text{\,if\quad}3<n \\
2tx_{12}x_{13}+t^2x_{12}^2,&\text{\,if\quad}n=3.
\end{array}
\right.
\end{equation}
For $(x_{2n}+x_{2n}^{-1})^{E_{23}(t)}$ we have
$$
(x_{2n}+x_{2n}^{-1})^{E_{23}(t)}=-\left(\sum_{r=3}^{n-1}x_{2r}x_{rn}^{-1}\right)
^{E_{23}(t)}=-\left(x_{23}x_{3n}^{-1}+\sum_{r=4}^{n-1}x_{2r}x_{rn}^{-1}\right)^{E_{23}(t)}=
$$
$$
-\left([x_{23}+t]x_{3n}^{-1}+\sum_{r=4}^{n-1}x_{2r}x_{rn}^{-1}\right)=
x_{2n}+x_{2n}^{-1}-tx_{3n}^{-1}.
$$
Since
$$
(x_{2n}- x_{2n}^{-1})^{E_{23}(t)}=[2x_{2n}-(x_{2n}+
x_{2n}^{-1})]^{E_{23}(t)}=[2x_{2n}-(x_{2n}+
x_{2n}^{-1}-tx_{3n}^{-1})]
$$
$$
=x_{2n}-x_{2n}^{-1}+tx_{3n}^{-1}
$$
we conclude that
\begin{equation}
\label{} (x_{2n}\pm x_{2n}^{-1})^{E_{23}(t)}=x_{2n}\pm x_{2n}^{-1}
\mp tx_{3n}^{-1}.
\end{equation}
Finally we have
$$
w_{2n}(xE_{23}(t))-w_{2n}(x)=(x_{2n}+x_{2n}^{-1}-tx_{3n}^{-1})
(x_{2n}-x_{2n}^{-1}+tx_{3n}^{-1})-(x_{2n}+x_{2n}^{-1})(x_{2n}-x_{2n}^{-1})=
$$
$$
tx_{3n}^{-1}(x_{2n}+x_{2n}^{-1}+x_{2n}-x_{2n}^{-1})-t^2(x_{3n}^{-1})^2=
2tx_{2n}^{-1}x_{3n}^{-1}-t^2(x_{3n}^{-1})^2,
$$
\begin{equation}
\label{w_2n^23}
 w_{2n}(xE_{23}(t))-w_{2n}(x)=
2tx_{2n}^{-1}x_{3n}^{-1}-t^2(x_{3n}^{-1})^2.
\end{equation}
Using (\ref{w_1n^23}) and (\ref{w_2n^23}) we get
\begin{equation}
\label{w_kn^23} w_{kn}(xE_{23}(t))-w_{kn}(x)=
\left\{\begin{array}{ll}
2tx_{12}x_{13}+t^2x_{12}^2,                 &\text{\,if\quad}n=3,\,\,k=1\\
2tx_{2n}^{-1}x_{3n}^{-1}-t^2(x_{3n}^{-1})^2,&\text{\,if\quad}k=2,\,\,n\geq 4 \\
0,                                          &\text{\,otherwise}.
\end{array}
\right.
\end{equation}
At last using (\ref{[T_23,De]})  and (\ref{w_kn^23}) we have
$$
\{T_{23}(t),\Delta^{is}(x)\}=\exp \left( -is\left[
b_{13}(2tx_{12}x_{13}+t^2x_{12}^2)+ \sum_{n=4}^\infty
b_{2n}(2tx_{2n}^{-1}x_{3n}^{-1}-t^2(x_{3n}^{-1})^2)
\right]\right).
$$
Further we get
\begin{equation}
\label{} \{T_{13}(t_2)\{T_{23}(t_1),\Delta^{is}(x)\}\}=
\exp\left(-isb_{13}2t_1t_2x_{12}\right).
\end{equation}
Indeed
$$
\{T_{13}(t_2)\{T_{23}(t_1),\Delta^{is}(x)\}\}=
$$
$$
 \exp\left(-is
b_{13}\left[(2t_1x_{12}[x_{13}+t_2]-t_1^2x_{12}^2)-
(2t_1x_{12}x_{13}-t_1^2x_{12}^2)\right]\right)
$$
$$
= \exp\left(-isb_{13}2t_1t_2x_{12}\right),
$$
compare with (\ref{[13[23,D]]}): $-[A^R_{13},[A^R_{23}
,\,\ln\Delta]]=2b_{13}x_{12}$! We have $\exp(itx_{12})\in M_\phi'$
and hence $\exp(itx_{12}^2)\in M_\phi'$. Using expression for
$\{T_{23}(t_1),\Delta^{is}(x)\}$ we conclude that
$$
M_\phi'\ni\{T_{23}(t_1),\Delta^{is}(x)\}\exp(isb_{13}t^2x_{12}^2)=
$$
$$
\exp \left( -is\left[ b_{13}(2tx_{12}x_{13})+ \sum_{n=4}^\infty
b_{2n}(2tx_{2n}^{-1}x_{3n}^{-1}-t^2(x_{3n}^{-1})^2)
\right]\right),
$$
so
$$
M_\phi'\ni\{T_{12}(t_2),\{T_{23}(t_1),\Delta^{is}(x)\}\exp(isb_{13}t^2x_{12}^2)\}
= \exp\left(-isb_{13}2t_1t_2x_{13}\right).
$$
Compare with the expression
$-[A^R_{12},[A^R_{23},\,\ln\Delta]]=2b_{13}x_{13} $. Finally we
conclude that
\begin{equation}
\label{12,13in M'} \exp(itx_{12}),\quad \exp(itx_{13})\in M_\phi'
\end{equation}
In general (without any projections) the following lemma holds
\begin{lem}  We have
\begin{equation}
\label{{T^R(t),w}}
w_{kn}(xE_{mm+1}(t))-w_{kn}(x)=
\left\{\begin{array}{ll}
2tx_{rm}x_{rm+1}+t^2x_{rm+1}^2,&\text{\,if\,}\,\,n=m+1,\,\,1\leq k\leq m\!-\!1\\
2tx_{mn}^{-1}x_{m+1n}^{-1}-t^2(x_{m+1n}^{-1})^2,&\text{\,if\,}\,\,k=m,\,\,n\geq m+2\\
0,                        &\text{\,otherwise},
\end{array}\right.
\end{equation}
hence
$$
\{T_{mm+1}(t),\Delta^{is}(x)\}=
$$
\begin{equation}
\label{[T_kn,De]} \exp\left(-is\left[
\sum_{r=1}^{m-1}b_{rm+1}(2tx_{rm}x_{rm+1}+t^2x_{rm+1}^2)+
\sum_{n=m+2}^\infty b_{mn}
(2tx_{mn}^{-1}x_{m+1n}^{-1}-t^2(x_{m+1n}^{-1})^2)
 \right] \right).
\end{equation}
\end{lem}
\begin{pf} The proof is similar to the proof of the Lemma \ref{l.[A^R,w]}.
\qed\end{pf}
To obtain another functions $\exp(itx_{kn})$ in the general case
we should make all the steps as it was indicated before. For
example to obtain $\exp(is x_{14}),\,\,\exp
(isx_{23}),\,\,\exp(isx_{24})$ we should do  {\bf the second step}
i.e. consider the operators
$$
\{T_{34}(t),\Delta^{is}(x)\}
$$
and all necessary combinations.

To obtain $\exp(is x_{15}),\,\,\exp
(isx_{25}),\,\,\exp(isx_{34}),\,\exp(isx_{34})$ we should consider
the following operators
$$
\{T_{45}(t),\Delta^{is}(x)\},
$$
and so on. Finally we shall obtain all functions
$\exp(isx_{kn}),\,k<n$.

\section{ Example of the measure} We show that the set $b=(b_{kn})_{k<n}$
 for which
$$
S^{L}_{kn}(b)<\infty,\quad E(b)<\infty,\quad
\text{\,and\,}\,\,S^{R,L}_{kn}(b)=\infty,\,\,\,\,1\leq k<n,
$$
 where
$$
S^{L}_{kn}(b)=\sum_{m=n+1}^\infty \frac{b_{km}}{b_{nm}},\,\,
E(b)=\sum_{k<n}\frac{S^{L}_{kn}(b)}{b_{kn}},\,\,
S^{R,L}_{kn}(b)=\sum_{m=n+1}^\infty \frac{b_{km}}{S^{L}_{nm}(b)}.
$$
is not empty. Indeed let us take $b_{kn}=(a_k)^n.$ We have
$$
S^{L}_{kn}(b)=\sum_{m=n+1}^\infty \left(\frac{a_k}{a_n}\right)^m=
\left(\frac{a_k}{a_n}\right)^{n+1}\sum_{m=0}^\infty
\left(\frac{a_k}{a_n}\right)^m=\left(\frac{a_k}{a_n}\right)^{n+1}\frac{1}{1-\frac{a_k}{a_n}}<\infty
$$
iff $a_k<a_{k+1},\,\,k\in {\mathbb N}$, for example $a_k=s^k$ with
$s>1$. Further we get
$$
E(b)=\sum_{k=1}^\infty\sum_{n=k+1}^\infty\frac{S^{L}_{kn}(b)}{b_{kn}}=
\sum_{k=1}^\infty\sum_{n=k+1}^\infty
\left(\frac{a_k}{a_n}\right)^{n+1}\frac{1}{1-\frac{a_k}{a_n}}\frac{1}{a_k^n}=
$$
$$
\sum_{k=1}^\infty a_k \sum_{n=k+1}^\infty
\left(\frac{1}{a_n}\right)^{n+1}\frac{1}{1-\frac{a_k}{a_n}}
<\sum_{k=1}^\infty
\frac{a_k}{1-\frac{a_k}{a_{k+1}}}\sum_{n=k+1}^\infty
\left(\frac{1}{a_n}\right)^{n+1}
$$
$$
<\sum_{k=1}^\infty
\frac{a_k}{1-\frac{a_k}{a_{k+1}}}\sum_{n=k+1}^\infty
\left(\frac{1}{a_{k+1}}\right)^{n+1}
=\sum_{k=1}^\infty\frac{a_k}{1-\frac{a_k}{a_{k+1}}}
\left(\frac{1}{a_{k+1}}\right)^{k+2}
\frac{1}{1-\frac{1}{a_{k+1}}}=
$$
$$
\sum_{k=1}^\infty\frac{\frac{a_k}{a_{k+1}}}{1-\frac{a_k}{a_{k+1}}}
\left(\frac{1}{a_{k+1}}\right)^{k} \frac{1}{a_{k+1}-1}
 <
\sum_{k=1}^\infty\frac{\frac{a_k}{a_{k+1}}}{1-\frac{a_k}{a_{k+1}}}
\left(\frac{1}{a_{2}}\right)^{k} \frac{1}{a_{2}-1}.
$$
If for example $a_k=s^k$ with $s>1$ we have
$$
E(b)<\frac{\frac{1}{s}}{1-\frac{1}{s}}\sum_{k=1}^\infty
\frac{1}{s^{k(k+1)}}
\frac{1}{s^{k+1}-1}<\infty.
$$
At last
$$
S^{R,L}_{kn}(b)=\sum_{m=n+1}^\infty \frac{b_{km}}{S^{L}_{nm}(b)}=
\sum_{m=n+1}^\infty \frac{a_k^m\left(1-\frac{a_n}{a_m}\right)}{
\left(\frac{a_n}{a_m}\right)^{m+1} }
$$
$$
=\sum_{m=n+1}^\infty
\left(\frac{a_ka_m}{a_n}\right)^m\left(\frac{a_m}{a_n}\right)
\left(1-\frac{a_n}{a_m}\right)=\sum_{m=n+1}^\infty
\left(\frac{a_ka_m}{a_n}\right)^m
\left(\frac{a_m}{a_n}-1\right)=\infty,
$$
if $\lim_ma_m=\infty.$ For $a_k=s^k$ with $s>1$ we have
$$
S^{R,L}_{kn}(b)=\sum_{m=n+1}^\infty
s^{(m+k-n)m}(s^{m-n}-1)=\infty.
$$
\section{Modular operator}
We recall how to find the modular operator and the operator of
canonical conjugation for the von Neumann algebra ${\mathfrak
A}^\rho_G$, generated by the right regular representation $\rho$
of a locally compact Lie group $G$. Let $h$ be a right invariant
Haar measure on $G$ and
$$
\rho,\lambda:G\mapsto U(L^2(G,h))
$$
be the right and the left regular representations of the group $G$
defined by
$$
(\rho_tf)(x)=f(xt),\,\,(\lambda_tf)(x)=(dh(t^{-1}x)/dh(x))^{-1/2}f(t^{-1}x).
$$
To define  the right Hilbert algebra on $G$ we can proceed as
follows. Let $M(G)$ be algebra  of all probability measures on $G$
with convolution
$$
(\mu*\nu)(s)=
$$
We define the homomorphism
$$
M(G)\ni\mu\mapsto \rho^\mu=\int_G\rho_t\d\mu(t)\in B(L^2(G,h)).
$$
We have $\rho^{\mu}\rho^{\nu}=\rho^{\mu*\nu}$, indeed
$$
\rho^{\mu}\rho^{\nu}=\int_G\rho_t\d\mu(t)\int_G\rho_s\d\nu(s)=
\int_G\int_G\rho_{ts}\d\mu(t)\nu(s)=\int_G\rho_t\d(\mu*\nu)(t)=\rho^{\mu*\nu}.
$$
Let us consider a subalgebra $M_h(G):=(\nu\in M(G)\mid \nu\sim h)$
of the algebra $M_h(G)$  In the case when $\mu\in M_h(G)$ we can
associate with the measure $\mu$ its Rodon-Nikodim derivative
$d\nu(t)/dh(t)=f(t)$. When $f\in C_0^\infty(G)$ or $f\in L^1(G)$
we can write
$$
\rho^f=\int_Gf(t)\rho_tdh(t),
$$
hence we can replace the algebra $M_h(G)$ by its subalgebra
identified with algebra of functions $C_0^\infty(G)$ or $L^1(G,h)$
with convolutions.
If we replace the Haar measure $h$ with some  measure $\mu\in
M_h(G)$ we obtain the isomorphic image $T^{R,\mu}$ of the right
regular representation $\rho$ in the space $L^2(G,\mu)$:
$T^{R,\mu}_t=U\rho_tU^{-1}$ where $U:L^2(G,h)\mapsto L^2(G,\mu)$
defined by $(Uf)(x)=\left(\frac{dh(x)}{d\mu(x)}\right)^{1/2}f(x)$.
we have
$$(T^{R,\mu}_tf)(x)=\left(\frac{d\mu(xt)}{d\mu(x)}\right)^{1/2}f(xt),$$
and
$$
T^f=\int_Gf(t)T^{R,\mu}_td\mu(t).
$$
We have (see \cite{{Con94}}, p.462) (we shall write $T_t$ instead
of $T^{R,\mu}_t$ )
$$
S(T^f):=(T^f)^*=\int_G\overline{f(t)}T_{t^{-1}}d\mu(t)=
\int_G\overline{f(t)}T_{t^{-1}}\frac{d\mu(t)}{d\mu(t^{-1})}d\mu(t^{-1})
$$
$$
\int_G\frac{d\mu(t^{-1})}{d\mu(t)}\overline{f(t^{-1})}T_{t}d\mu(t).
$$
Hence
$$
(Sf)(t)=\frac{d\mu(t^{-1})}{d\mu(t)}\overline{f(t^{-1})}.
$$
To calculate $S^*$ we use the fact that $S$ is antilinear so $
(Sf,g)=(S^*g,f)$. We have
$$
(Sf,g)=\int_G\frac{d\mu(t^{-1})}{d\mu(t)}\overline{f(t^{-1})}\overline{g(t)}d\mu(t)=
\int_G\overline{f(t^{-1})}\overline{g(t)}d\mu(t^{-1})=
$$
$$
\int_G\overline{g(t^{-1})}\overline{f(t)}d\mu(t)=(S^*g,f),
$$
hence $(S^*g)(t)=\overline{g(t^{-1})}.$ Finally the modular
operator $\Delta$ defined by $\Delta=S^*S$ has the following form
$(\Delta f)(t)=\frac{d\mu(t)}{d\mu(t^{-1})}f(t)$. Indeed we have
$$
f(t) \stackrel{S}{\mapsto}
\frac{d\mu(t^{-1})}{d\mu(t)}\overline{f(t^{-1})}\stackrel{S^*}{\mapsto}
\frac{d\mu(t)}{d\mu(t^{-1})}f(t).
$$
Finally , since $J=S\Delta^{-1/2}$ (see \cite{{Con94}} p.462) we
get
$$
f(t)
\stackrel{\Delta^{-1/2}}{\mapsto}\left(\frac{d\mu(t^{-1})}{d\mu(t)}\right)^{1/2}f(t)
\stackrel{J}{\mapsto}\frac{d\mu(t^{-1})}{d\mu(t)}
\left(\frac{d\mu(t)}{d\mu(t^{-1})}\right)^{1/2}\overline{f(t^{-1})}
$$
$$
=
\left(\frac{d\mu(t^{-1})}{d\mu(t)}\right)^{1/2}\overline{f(t^{-1})}.
$$
Hence
$$
(Jf)(t)=\left(\frac{d\mu(t^{-1})}{d\mu(t)}\right)^{1/2}\overline{f(t^{-1})},\,\,\text{\,and\,}\,\,
(\Delta f)(t)=\frac{d\mu(t)}{d\mu(t^{-1})}f(t).
$$
To prove that $JT^{R,\mu}_tJ=T^{L,\mu}_t$ we get
$$
f(t)\stackrel{J}{\mapsto}\left(\frac{d\mu(x^{-1})}{d\mu(x)}\right)^{1/2}\overline{f(x^{-1})}
\stackrel{T^{R,\mu}_t}{\mapsto}
\left(\frac{d\mu(xt)}{d\mu(x)}\right)^{1/2}
\left(\frac{d\mu((xt)^{-1})}{d\mu(xt)}\right)^{1/2}\overline{f((xt)^{-1})}=
$$
$$
\left( \frac{d\mu(t^{-1}x^{-1})}{d\mu(x)}
\right)^{1/2}\overline{f(t^{-1}x^{-1})}\stackrel{J}{\mapsto}
\left( \frac{d\mu(x^{-1})}{d\mu(x)} \right)^{1/2}\left(
\frac{d\mu(t^{-1}x)}{d\mu(x^{-1})} \right)^{1/2}f(t^{-1}x)=
$$
$$
\left(\frac{d\mu(t^{-1}x)}{d\mu(x)} \right)^{1/2} f(t^{-1}x)=
(T^{L,\mu}_tf)(x).
$$
\begin{rem} The representation $T^{R,\mu_b}$ is the inductive
limit of the representations $T^{R,\mu_b^m}$ of the group
$B(m,{\mathbb R})$ where the measure $\mu_b^m$ is the projection
of the measure $\mu_b$ onto subgroup   $B(m,{\mathbb R})$.
Obviously $\mu_b^m$ is equivalent with the Haar measure $h_m$ on
$B(m,{\mathbb R})$.
\end{rem}

\section{The uniqueness of the constructed factor}
Let $G$ be a solvable separable locally compact group or a
connected locally compact group. Then any representation $\pi$ of
$G$ in a Hilbert space generates an approximately
finite-dimensional von Neumann algebra (see \cite{Con76}).

Theorem 15 from V.9 p. 504 \cite{Con94} (Haagerup) There exists up
to isomorphism only one amenable factor of type $III_1$, the
factor $R_\infty$ of Araki and Woods (see \cite{Haa87}).

{\it Acknowledgements.} {The author would like to thank the
Max-Planck-Institute of Mathematics where the main part of the
work was done for the hospitality. He is grateful  to Prof.
A.Connes for  useful remarks concerning the uniqueness of the
constructed factor.
The partial financial support by the DFG project 436 UKR 113/87
 are gratefully acknowledged.}

\newpage

\end{document}